\documentclass{article}
\usepackage{graphicx} 
\usepackage{framed}
\usepackage[normalem]{ulem} 
\usepackage{enumerate}
\usepackage{tikz}
\usepackage{amssymb, amsmath, amsthm}
\usepackage{mathrsfs}
\usepackage{stmaryrd}
\usepackage[makeroom]{cancel}
\usepackage{graphicx,epstopdf}
\usepackage{gensymb}
\usepackage{array, tabu}
\usepackage[portrait, margin=1in]{geometry}
\usepackage{subcaption}
\usepackage{parskip}
\usepackage[english]{babel}  
\usepackage[backend=biber]{biblatex}  
\usepackage{authblk} 
\usepackage{hyperref}

\DefineBibliographyStrings{english}{
    andothers = {et al.},
    in        = {in },
    pages     = {pp.},
}

\DeclareFieldFormat{postnote}{#1} 
\DeclareNameAlias{sortname}{last-first} 
\DeclareFieldFormat[article, incollection]{title}{#1} 

\newcommand{\R}{\mathbb{R}}

\title{The Optimal Ratio of a Generalized Chaos Game in Regular Polytopes}
\author{Christoffer Tarmet}
\affil{Mentor: Dr. Markus Reitenbach}
\affil{Department of Mathematics, Colorado Mesa University}
\date{December 2024}

\begin{document}

\maketitle

\section{Abstract}
This paper investigates the concept of an optimal ratio for regular polytopes in $n$-dimensional space within the framework of the Generalized Chaos Game. The optimal ratio, $r_{\text{opt}}$, is defined as the value at which the self-similar regions of the resulting fractal touch but do not overlap. Using a series of Python simulations, we explore how the optimal ratio varies across different polytopes, from two-dimensional polygons to three-dimensional polyhedra and beyond. The results, visualized through plots generated for various polytopes and values of the scaling factor $r$, demonstrate that the optimal ratio is not universal but rather depends on each polytope's specific properties. A formula is then derived for determining the optimal ratio for any regular polytope in any dimension. The formula is then experimentally verified using multiple Python programs designed to search and find the optimal ratio iteratively. 

\newpage 

\section{Introduction} 
Fractals have captivated mathematicians, scientists, and artists alike for their self-similar structures and often intricate, seemingly infinite detail. Fractals have come to represent a class of geometric objects defined not by smooth boundaries, but by recursive, repeating patterns that remain complex under magnification. These structures appear throughout nature, from the branching of trees to the contours of coastlines, and have led to numerous applications in fields as diverse as biology, physics, and computer graphics.

One notable method for generating fractals is the Chaos Game, a stochastic process in which points are plotted iteratively, converging to reveal fractal structures within a given shape. Originally applied to two-dimensional polygons, the Chaos Game has traditionally yielded classic fractals such as the Sierpinski Triangle. In this paper, we extend the Chaos Game into higher dimensions, investigating its behavior within regular polytopes, including those beyond three-dimensional polyhedra. By defining and exploring an optimal ratio—a specific scaling factor where self-similar regions of the fractal touch without overlapping—our work offers insights into the conditions that generate distinct fractal patterns in various polytopes.

The goal of this study is twofold: to generalize the Chaos Game to higher-dimensional regular polytopes and to develop a conjecture that determines the optimal ratio for any regular polytope. Through a combination of theoretical analysis and Python-based simulations, we aim to offer a framework for understanding fractals in higher-dimensional spaces, shedding light on both their underlying mathematics and their potential applications.
\section{Chaos Game}
Generally, a chaos game is described as a method to generate fractals through an iterative random process. The rules of such a process are as follows:
\begin{enumerate}
    \item Pick a random point $c$ inside a regular $n$-gon and plot it. Also, pick a ratio $r\in(0.1)$. This ratio will remain fixed.
\item Select a random vertex of the $n$-gon.
\item Calculate the vector $v$ from the current point $c$ to the randomly selected vertex. 
\item Plot the next point $c+rv$.
\item Set the point $c + rv$ as the current point and repeat from step 2 \cite{rules}.
\end{enumerate}
Once the select number of iterations is completed, the initial few points are discarded. This process is known to generate fractals depending on the polygon and the value of $r$ \cite{rules}. 

\subsection{Generalized Chaos Game}
In this paper, we will modify the rules above in order to include higher-dimensional polytopes. Let $P_n$ be an $n$-dimensional regular polytope with vertices $v_1,v_2, \ldots, v_V$. Then, the Generalized Chaos Game (GCG) can be described as generating a sequence of points in $\R^n$ using the iterative function: \begin{align}
    x_{c+1} = (1-r)x_c + rv_i
\end{align}
where $x_0$ is a randomly selected initial point, and $v_i$ is one of the polytope's vertices with $i\in \{1,2,\ldots, V\}$, chosen at random in each iteration. To visually represent the Chaos Game, we plot the points generated by this iteration. Each point is assigned a color based on which vertex $v_i$ was used in that iteration. The reader should note that the iterative function (1), given a 2-dimensional polytope (i.e. a polygon) is equivalent to the procedure discussed in section 3. However, the power of the GCG as defined is that it can play in any $n$-dimensional regular polytope.

\subsection{Plots of GCG in two and three-dimensional polygons}
The following plots were generated from a GCG played inside two-dimensional polygons and three-dimensional polyhedra, using different values of $r$. These plots were created using the Python programs provided in \cite{programs} where the initial 6 points were discarded.

\begin{figure}[h!]
\centering
\begin{subfigure}{0.3\textwidth}
  \centering
  \includegraphics[width=1\linewidth]{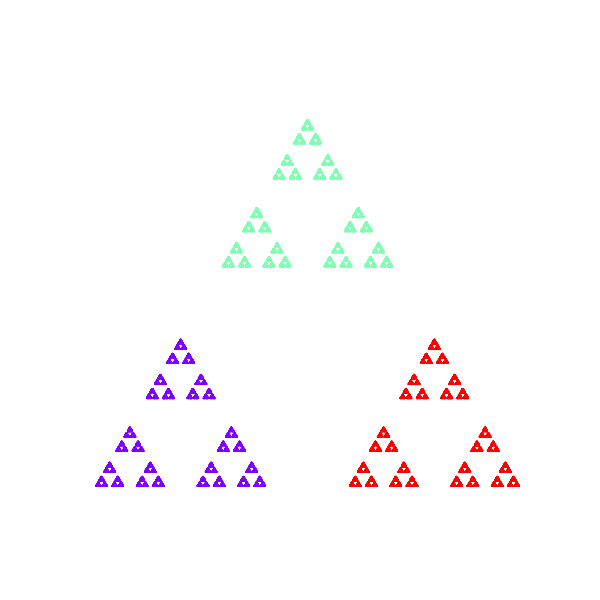}
  \caption{$r = 0.6$}
  \label{fig:sub1}
\end{subfigure}%
\begin{subfigure}{0.3\textwidth}
  \centering
  \includegraphics[width=1\linewidth]{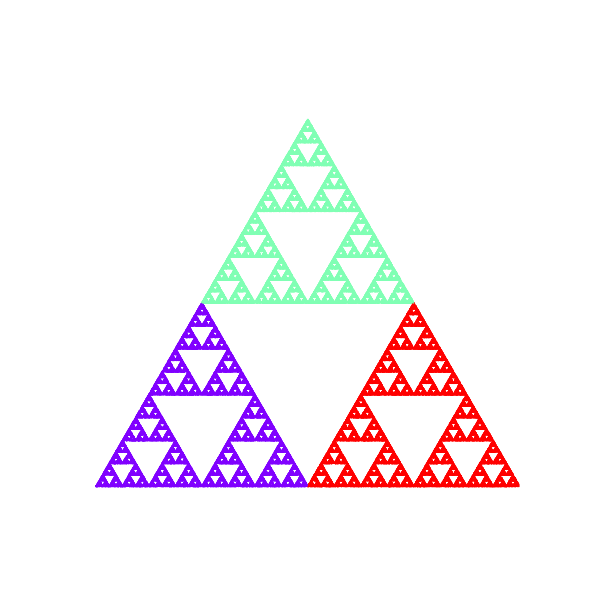}
  \caption{$r = 0.5$}
  \label{fig:sub2}
\end{subfigure}
\begin{subfigure}{0.3\textwidth}
    \centering
  \includegraphics[width=1\linewidth]{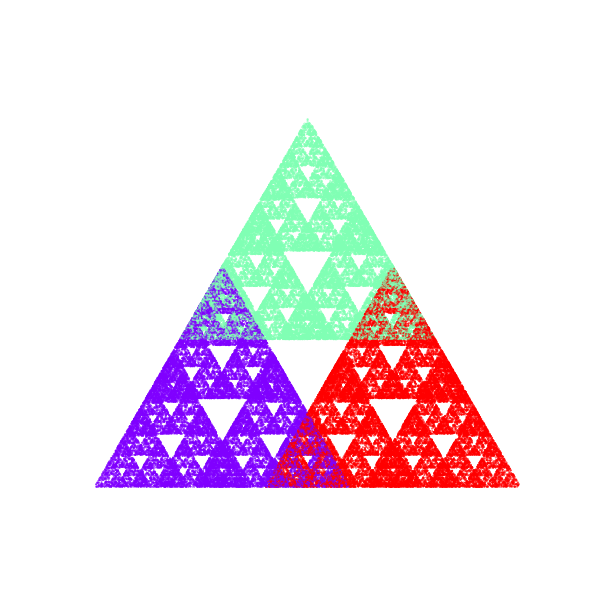}
  \caption{$r = 0.4$}
  \label{fig:sub2}
\end{subfigure}
\label{fig:test}
\caption{Equilateral triangles with different values of $r$}
\end{figure}

\begin{figure}[h!]
\centering
\begin{subfigure}{0.3\textwidth}
  \centering
  \includegraphics[width=1\linewidth]{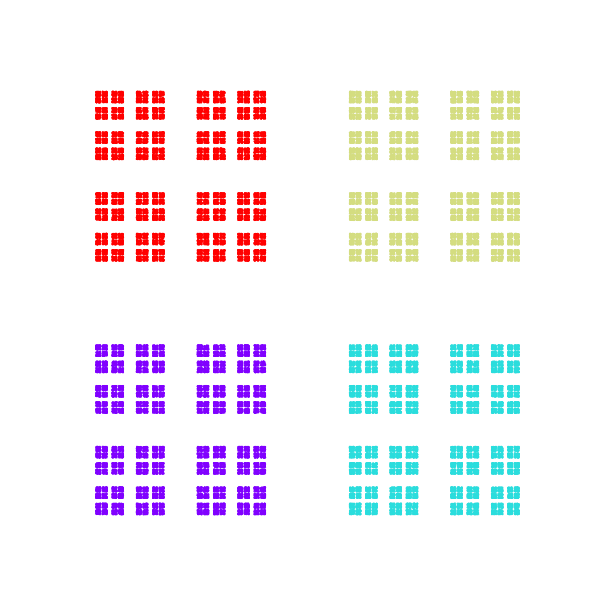}
  \caption{$r = 0.6$}
  \label{fig:sub1}
\end{subfigure}%
\begin{subfigure}{0.3\textwidth}
  \centering
  \includegraphics[width=1\linewidth]{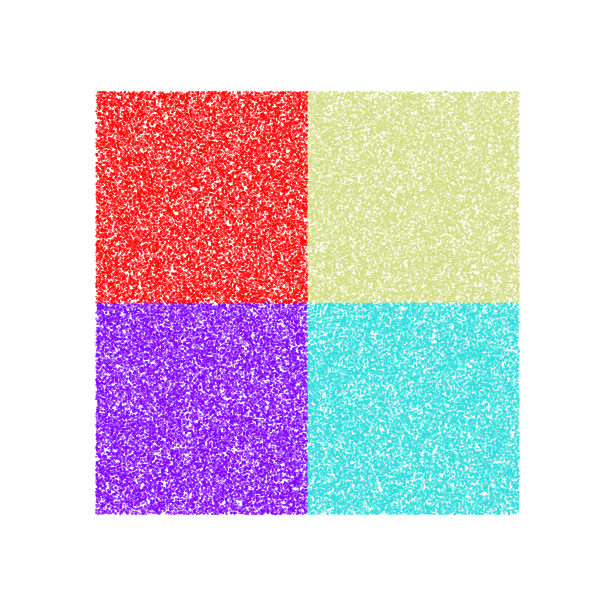}
  \caption{$r = 0.5$}
  \label{fig:sub2}
\end{subfigure}
\begin{subfigure}{0.3\textwidth}
    \centering
  \includegraphics[width=1\linewidth]{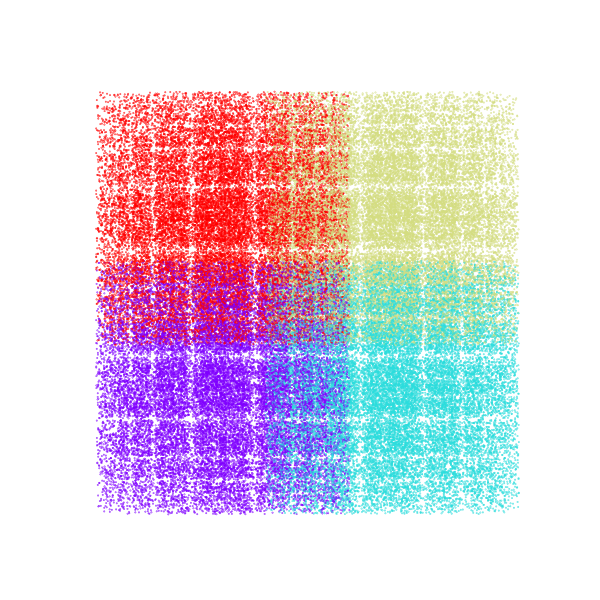}
  \caption{$r = 0.4$}
  \label{fig:sub2}
\end{subfigure}
\label{fig:test}
\caption{Squares with different values of $r$}
\end{figure}

\begin{figure}[h!]
\centering
\begin{subfigure}{0.3\textwidth}
  \centering
  \includegraphics[width=1\linewidth]{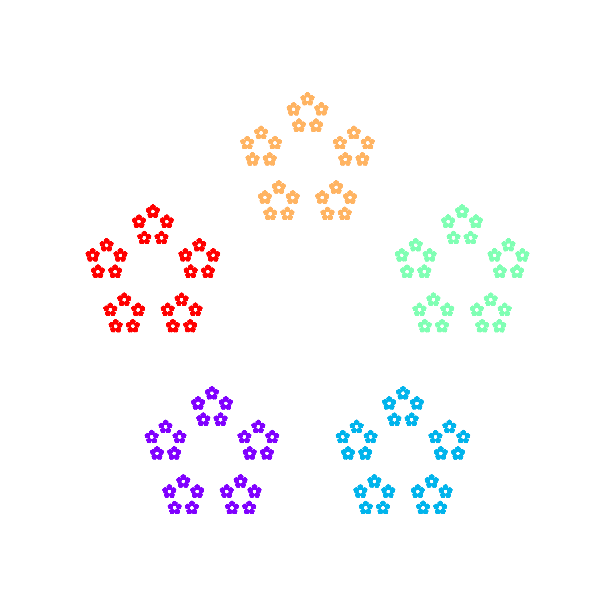}
  \caption{$r = 0.7$}
  \label{fig:sub1}
\end{subfigure}%
\begin{subfigure}{0.3\textwidth}
  \centering
  \includegraphics[width=1\linewidth]{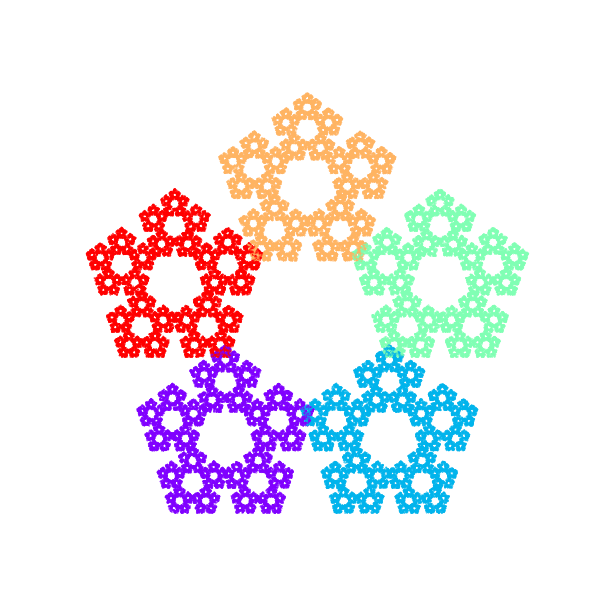}
  \caption{$r = 0.6$}
  \label{fig:sub2}
\end{subfigure}
\begin{subfigure}{0.3\textwidth}
    \centering
  \includegraphics[width=1\linewidth]{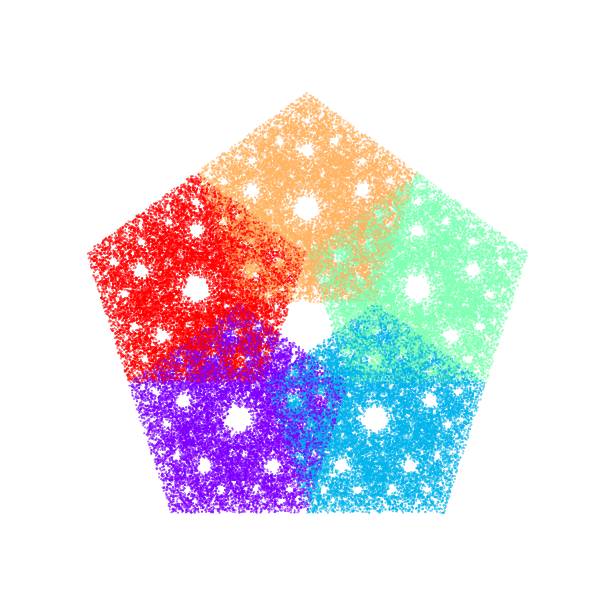}
  \caption{$r = 0.5$}
  \label{fig:sub2}
\end{subfigure}
\label{fig:test}
\caption{Regular pentagons with different values of $r$}
\end{figure}

\begin{figure}[h!]
\centering
\begin{subfigure}{0.3\textwidth}
  \centering
  \includegraphics[width=1\linewidth]{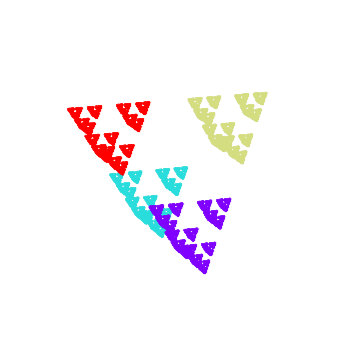}
  \caption{$r = 0.6$}
  \label{fig:sub1}
\end{subfigure}%
\begin{subfigure}{0.3\textwidth}
  \centering
  \includegraphics[width=1\linewidth]{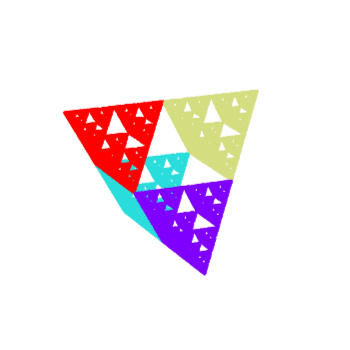}
  \caption{$r = 0.5$}
  \label{fig:sub2}
\end{subfigure}
\begin{subfigure}{0.3\textwidth}
    \centering
  \includegraphics[width=1\linewidth]{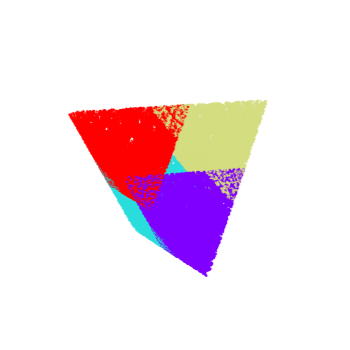}
  \caption{$r = 0.4$}
  \label{fig:sub2}
\end{subfigure}
\label{fig:test}
\caption{Regular tetrahedra with different values of $r$}
\end{figure}

\begin{figure}[h!]
\centering
\begin{subfigure}{0.3\textwidth}
  \centering
  \includegraphics[width=1\linewidth]{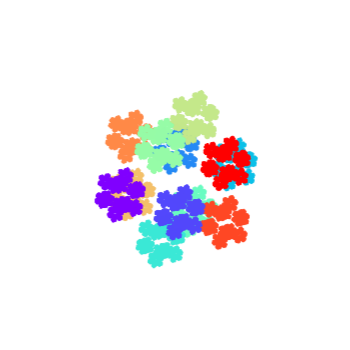}
  \caption{$r = 0.7$}
  \label{fig:sub1}
\end{subfigure}%
\begin{subfigure}{0.3\textwidth}
  \centering
  \includegraphics[width=1\linewidth]{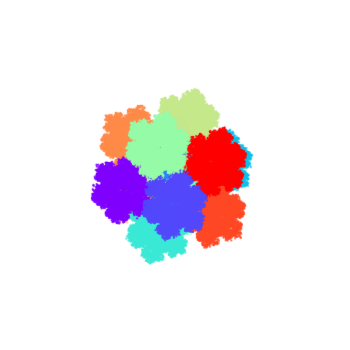}
  \caption{$r = 0.6$}
  \label{fig:sub2}
\end{subfigure}
\begin{subfigure}{0.3\textwidth}
    \centering
  \includegraphics[width=1\linewidth]{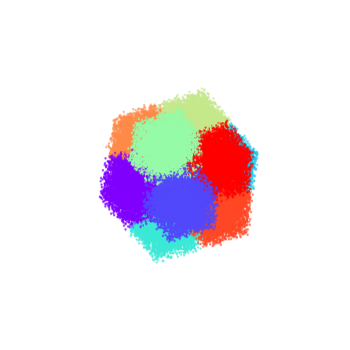}
  \caption{$r = 0.5$}
  \label{fig:sub2}
\end{subfigure}
\label{fig:test}
\caption{Regular icosahedra with different values of $r$}
\end{figure}

\textbf{Definition 1:} The optimal ratio of a polytope, $r_{\text{opt}}$, is the ratio at which the self-similar regions of the fractal touch without overlapping. In other words, these regions may share vertex coordinates, but no point from one region appears within another region of a different color.

\vspace{10pt}

This paper will examine the optimal ratio for regular polytopes in 
$n$ dimensions. As shown in Figures 1 through 5, the optimal ratio can vary depending on the specific polytope.

\newpage

\section{Deriving the Optimal Ratio Formula for a Generalized Chaos Game in Regular Polytopes}
To find a formula for the optimal ratio depending on the polytope, we start by describing the Chaos Game from the perspective of metric spaces. 

\vspace{10pt}
\textbf{Definition 2:}
Let $(\mathbf{X},d)$ be a complete metric space. Then $\mathcal{H}(\mathbf{X})$ denotes the space whose points are the compact subsets of $\mathbf{X}$, other than the empty set \cite{barnsley}.

\textbf{Definition 3:} 
Let $(\mathbf{X},d)$ be a complete metric space, $x\in\mathbf{X}$, and $B\in\mathcal{H}(\mathbf{X})$. Define \[d(x,B) = \text{min}\{d(x,y):y\in B\}.\]
$d(x,B)$ is called the distance from the point $x$ to the set $B$ \cite{barnsley}.

\textbf{Definition 4:} 
Let $(\mathbf{X},d)$ be a complete metric space. Let $A,B\in\mathcal{H}(\mathbf{X})$. Define \[d(A,B) = \text{max}\{d(x,B):x\in A\}.\]
$d(A,B)$ is called the distance from the set $A\in\mathcal{H}(\mathbf{X})$ to the set $B\in\mathcal{H}(\mathbf{X})$ \cite{barnsley}.

\textbf{Definition 5:} 
Let $(\mathbf{X}, d)$ be a complete metric space. Then the Hausdorff distance between points $A,B\in\mathcal{H}(\mathbf{X})$ is defined by \[h(A,B)=\max \{d(A,B), d(B,A)\}\]\cite{barnsley}.

\textbf{Definition 6:} 
    A transformation $w:\mathbf{X}\rightarrow \bold{X}$ on a metric space $(\mathbf{X}, d)$ is called
contractive or a contraction mapping if there is a constant $0 \leq  s < 1$ such that
\[d(w(x),w(y)) \leq s\cdot d(x,y) \hspace{10pt} \forall x,y\in \bold{X}\]
Any such number $s$ is called a contractivity factor for $w$ \cite{barnsley}.

\vspace{10pt}

For each vertex of the polytope, we can interpret a point being plotted towards that vertex as a contraction mapping. That is, if $P_n$ be a regular $n$-dimensional polytope with $V$ vertices denoted $v_1,v_2,\ldots, v_V\in\R^n$,
    let $(\bold{X},d)$ be the metric space with $\bold{X}=\R^n$ and $d$ is the Euclidean metric.
    Let $w_i:\bold{X}\rightarrow \bold{X}$ for each vertex $v_i$, where each $w_i$ maps points closer to $v_i$ via \[w_i(x) = (1-r)x + rv_i\]
    where $r\in(0,1)$ and $i\in\{1,2,\ldots, V\}$.
    
\vspace{10pt}

Note that each $w_i$ is a contraction mapping of $\bold{X}$ since \begin{align*}
    d(w_i(x),w_i(y)) &= \|w_i(x)-w_i(y)\| \\
    &= \|(1-r)x + rv_i -(1-r)y + rv_i\| \\
    &= (1-r)\cdot \|x-y\|.
\end{align*}
Since $(1-r)\in(0,1)$, each $w_i$ is a contraction mapping with contractivity factor $1-r$.

However, we are interested in the boundaries of these mappings. Therefore, we can define each contraction mapping as an affine transformation of the space defined by the boundary of the polytope to a scaled-down version where the only vertex coordinate that remains constant is $v_i$. 

\vspace{10pt}

\textbf{Lemma 1:} Let $w:\mathbf{X} \rightarrow \mathbf{X}$ be a contraction mapping on the metric space $(\mathbf{X}, d)$. Then $w$ is continuous \cite{barnsley}.

\vspace{10pt}

\textbf{Lemma 2:} Let $w:\mathbf{X} \rightarrow \mathbf{X}$ be a continuous mapping on the metric space $(\mathbf{X}, d)$ with contractivity factor $s$. Then $w: \mathcal{H}(\mathbf{X}) \rightarrow \mathcal{H}(\mathbf{X})$ defined by \[w(B) = \{w(x):x\in B\}\hspace{10pt}\forall B \in \mathcal{H}(\mathbf{X})\]
is a contraction mapping on $(\mathcal{H}(\mathbf{X}), h(d))$ with contractivity factor $s$, where $h(d)$ is the Hausdorff metric \cite{barnsley}. 

\vspace{10pt}

\textbf{Lemma 3:} Let $(\mathbf{X}$, d$)$ be a metric space. Let $\left\{w_i: i=1,2, \ldots, V\right\}$ be contraction mappings on $(\mathcal{H}(\mathbf{X}), h)$. Let the contractivity factor be $s$ for all $w_i$. Define $W: \mathcal{H}(\mathbf{X}) \rightarrow \mathcal{H}(\mathbf{X})$ by

\begin{align*}
W(B) & =w_1(B) \cup w_2(B) \cup \ldots \cup w_V(B) \\
& =\cup_{i=1}^V w_i(B), \quad \text { for each } B \in \mathcal{H}(\mathbf{X}).
\end{align*}

Then $W$ is a contraction mapping with contractivity factor $s$ \cite{barnsley}.

\vspace{10pt}
Let $W^n$ denote the $n$th iteration of applying the contraction mapping $W$. It follows from Lemma 3 that each iteration of $W$ will result in $V$ self-similar regions, one per vertex. $W^n$ can then be thought of as the possible space where a point can be plotted in iteration $n$ of the GCG. The first five iterations of applying $W$ to an equilateral triangle can be seen in Figure 6. 

\begin{figure}[h!]
\centering
\begin{subfigure}{0.3\textwidth}
  \centering
  \includegraphics[width=1\linewidth]{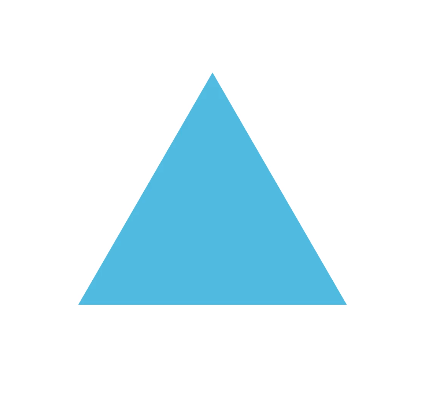}
  \caption{$W^0(B)$}
  \label{fig:sub1}
\end{subfigure}%
\begin{subfigure}{0.3\textwidth}
  \centering
  \includegraphics[width=1\linewidth]{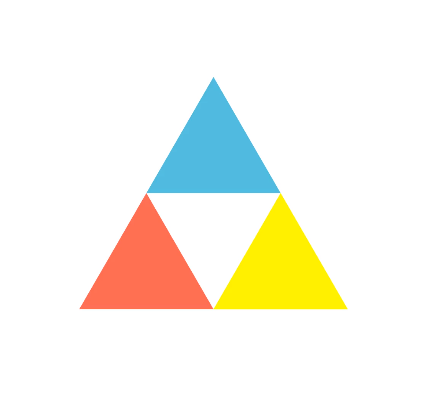}
  \caption{$W^1(B)$}
  \label{fig:sub2}
\end{subfigure}
\begin{subfigure}{0.3\textwidth}
    \centering
  \includegraphics[width=1\linewidth]{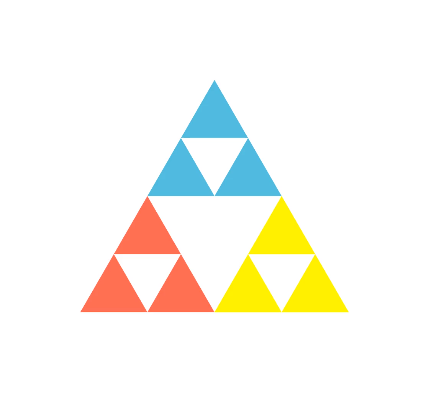}
  \caption{$W^2(B)$}
  \label{fig:sub2}
\end{subfigure}
\label{fig:test}
\begin{subfigure}{0.3\textwidth}
  \centering
  \includegraphics[width=1\linewidth]{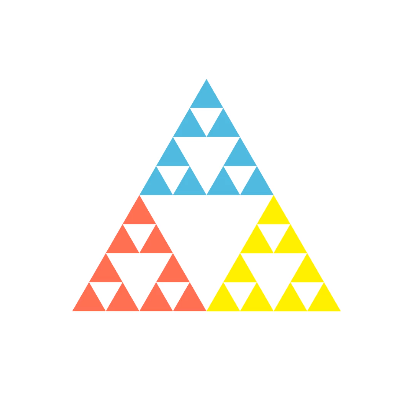}
  \caption{$W^3(B)$}
  \label{fig:sub1}
\end{subfigure}%
\begin{subfigure}{0.3\textwidth}
  \centering
  \includegraphics[width=1\linewidth]{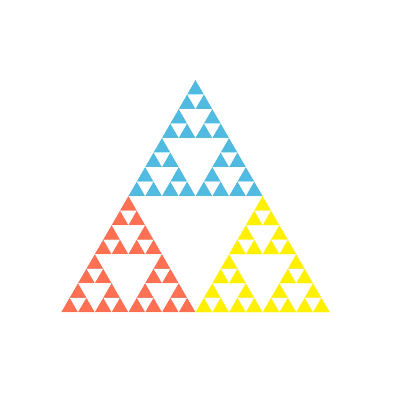}
  \caption{$W^4(B)$}
  \label{fig:sub2}
\end{subfigure}
\begin{subfigure}{0.3\textwidth}
    \centering
  \includegraphics[width=1\linewidth]{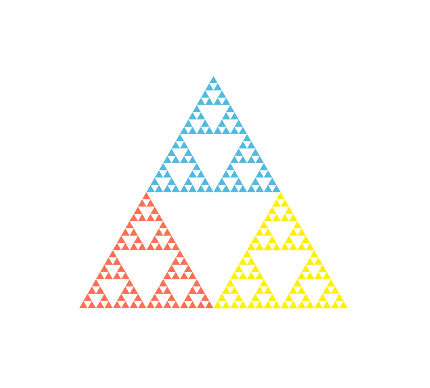}
  \caption{$W^5(B)$}
  \label{fig:sub2}
\end{subfigure}
\label{fig:test}
\caption{First Iterations of $W(B)$ on an Equilateral Triangle with $r=0.5$ \cite{programs}}
\end{figure}
Figure 6 can be thought of as three distinct affine transformations of the space where the space is scaled by a factor of $1-r$ and then translated to the vertices of the original equilateral triangle. 

We are interested in the optimal ratio of the Chaos Game. We want to determine the value of $r$ at which the differently colored self-similar regions, produced by applying $W$ to the regular $n$-dimensional polytope, touch each other without overlapping. Note that the boundaries of the self-similar regions remain constant through all iterations of $W$. Thus, the optimal ratio will also remain constant through all iterations. 

Consider the vectors connecting each pair of a polytope's vertices. There are 
    $V(V-1)$
possible vectors between these pairs, if $V$ is the number of vertices of the polytope. However, we are only interested in each unique vector direction, not both a vector and its negative. Thus the number of vectors decreases to $\binom{V}{2} = \frac{V(V-1)}{2}$. Figure 7 illustrates this concept with a regular pentagon in 2-dimensional space.

\begin{figure}[h!]
    \centering
    \includegraphics[width=0.32\linewidth]{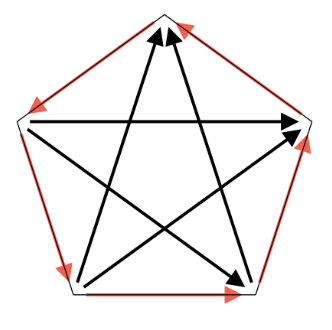}
    \caption{Visual representation of the vectors in $\mathcal{A}$ and $\mathcal{B}$ for a regular pentagon \cite{programs}}
    \label{fig:enter-label}
\end{figure}

We will use the notation $\Vec{u} \parallel \Vec{v}$ to denote two nonzero parallel vectors, $\Vec{u}, \Vec{v}$. That is, \begin{align*}
    \Vec{u}\parallel \Vec{v} \Leftrightarrow \Vec{u} = \lambda \Vec{v},
\end{align*} 
for some constant $\lambda$.

Let \( \mathcal{A} \) denote the set of all vectors between pairs of distinct vertices of the polytope, with each pair contributing exactly one vector (represented by the black and red vectors in Figure 7). For any pair of vertices, the vector included in \( \mathcal{A} \) is chosen to ensure that no two vectors in \( \mathcal{A} \) are opposites.

Let \( \mathcal{B} \) denote the subset of \( \mathcal{A} \) consisting of edge vectors (the red vectors in Figure 7). Note that the magnitudes of all vectors in \( \mathcal{B} \) are equal to the polytope's edge length, \( \ell \).

\textbf{Definition 8:} The maximum magnitude of one of the vectors in $\mathcal{A}$ that is parallel to at least one of the vectors in $\mathcal{B}$ is denoted $\delta_\parallel$ and is defined by
\begin{align*} 
    \delta_\parallel := \max\{|\vec{u}| : \vec{u} \parallel \vec{v}, \hspace{4pt} \vec{u} \in \mathcal{A} \hspace{4pt} \text{and} \hspace{4pt} \vec{v} \in \mathcal{B}\}.
\end{align*}
Note that $\delta_\parallel$ is well-defined since $\mathcal{B} \subseteq \mathcal{A}$.

\begin{figure}[h!]
    \centering
    \includegraphics[width=0.55\textwidth]{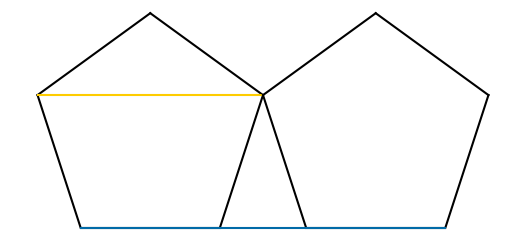}
    \caption{Two adjacent scaled-down polygons showing $(1-r_{\text{opt}})\delta_\parallel$ (yellow) and $\ell$ (blue) \cite{programs}}
    \label{fig:d_m}
\end{figure}

\newpage

\begin{figure} [h!]
    \centering
    \includegraphics[width=0.35\linewidth]{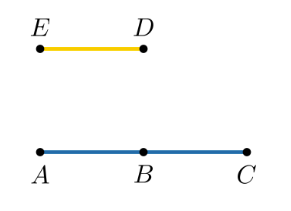}
    \caption{Sketch of the lines of length $(1-r_{\text{opt}})\delta_\parallel$ and $\ell$ if $\delta_\parallel = \ell$ \cite{programs} }
    \label{fig:enter-label}
\end{figure}

\underline{Case 1:} Suppose \( \delta_\parallel = \ell \), where \( \ell \) is the edge length of the polytope. In this case, the situation is depicted in Figure 9, where \( A \) and \( C \) are vertices of the original polytope, and \( B \) is the midpoint of \( \overline{AC} \). The points \( D \) and \( E \) represent the scaled-down vertices corresponding to the vertices that form the vector with the maximum magnitude in \( \mathcal{A} \), parallel to a vector in \( \mathcal{B} \), which defines \( \delta_\parallel \).

By the definition of \( \delta_\parallel \), the line segments $\overline{ED}$, with length \( |ED| = (1 - r_{\text{opt}}) \delta_\parallel \), and $\overline{AC}$, with length \( |AC| = \ell \), are parallel. The optimal ratio ensures that \( |ED| = |AB| \), where \( |AB| = \frac{1}{2} \ell \) since \( B \) is the midpoint of \( \overline{AC} \). Thus, we have the equation:
\[
(1 - r_{\text{opt}}) \delta_\parallel = \frac{1}{2} \ell.
\]
Solving for \( r_{\text{opt}} \), we get:
\[
r_{\text{opt}} = 1 - \frac{\ell}{2 \delta_\parallel}.
\]
Since \( \delta_\parallel = \ell \), this simplifies to:
\[
r_{\text{opt}} = \frac{1}{2}.
\]

\vspace{3cm}

\begin{figure}[h!]
    \centering
    \includegraphics[width=0.5\textwidth]{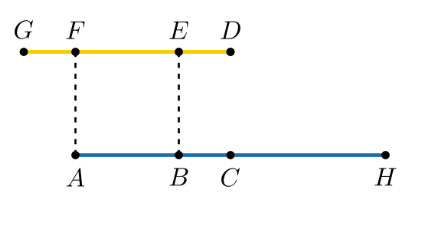}
    \caption{Sketch of the lines $(1-r_{\text{opt}})\delta_\parallel$ and $\ell$ if $\delta_\parallel \neq \ell$ \cite{programs}}
    \label{fig:d_m}
\end{figure}

\newpage

\underline{Case 2:} Suppose \( \delta_\parallel \neq \ell \). In this case, the situation is depicted in Figure 10, where \( A \) and \( H \) are vertices of the original polytope, and \( C \) is the midpoint of \( \overline{AH} \). Additionally, \( B \) is the scaled-down vertex corresponding to \( H \) when it is scaled and translated towards vertex \( A \).

The points \( D \) and \( E \) represent the scaled-down vertices corresponding to those that form the vector with the maximum magnitude in \( \mathcal{A} \), parallel to a vector in \( \mathcal{B} \), which defines \( \delta_\parallel \). The points \( F \) and \( E \) lie on \( \overline{GD} \) and are perpendicular to the line joining \( A \) and \( B \).

By the definition of \( \delta_\parallel \), the line segments $\overline{GD}$, with length \( |GD| = (1 - r_{\text{opt}}) \delta_\parallel \), and $\overline{AH}$, with length \( |AH| = \ell \), are parallel. From the figure, we have the following relationships:
\begin{itemize}
    \item \( |AB| = (1 - r_{\text{opt}}) \ell \)
    \item \( |FD| = |GD| - |GF| = |GD| - \frac{1}{2} (|GD| - |AB|) = \frac{1}{2}(1 - r_{\text{opt}})(\delta_\parallel + \ell) \)
\end{itemize}

The optimal ratio ensures that \( |FD| = |AC| \), so we equate:
\[
\frac{1}{2}(1 - r_{\text{opt}})(\delta_\parallel + \ell) = \frac{1}{2} \ell.
\]

Simplifying this expression and solving for \( r_{\text{opt}} \), we get:
\begin{align}
    r_{\text{opt}} = \frac{\delta_\parallel}{\delta_\parallel + \ell}.
\end{align}

Note that this formula for \( r_{\text{opt}} \) holds even when \( \delta_\parallel = \ell \), since in that case:
\[
r_{\text{opt}} = \frac{\ell}{\ell + \ell} = \frac{1}{2}.
\]

Thus, equation (2) provides a general formula for the optimal ratio of a Generalized Chaos Game (GCG) on any regular \( n \)-dimensional polytope.

\subsection{Finding $\delta_\parallel$ in Practice}
In practice, we find $\delta_\parallel$ by first orienting the polytope so that at least one of the edges is parallel to one of the $n$ axes. Without loss of generality, say the polytope is oriented such that one edge is parallel to the $x$-axis. In that case, $\delta_\parallel$ will be the largest distance between two vertices' $x$-coordinates. For example, say we have a regular icosahedron with edge length $\ell = 2$ with coordinates:
\begin{align*}
    (0, \pm \phi, \pm 1) \\
    (\pm 1, 0, \pm \phi) \\
    (\pm \phi, \pm 1, 0),
\end{align*}
where $\phi$ is the golden ratio \cite{MR151873}. 

Then, the polyhedron is oriented such that one edge is parallel to the $x$-axis since $(-1,0,\phi)$ and $(1,0,\phi)$ constitute an edge. Thus, $\delta_\parallel$ can be found by looking at the largest difference between any two vertices' $x$-coordinates. Hence, $\delta_\parallel = 2\phi$.

To find the optimal ratio for the icosahedron we use the formula derived above to obtain 
\begin{align*}
    r_{\text{opt}} = \frac{2\phi}{2\phi + 2} = \frac{\phi}{\phi + 1} = \frac{1}{\phi}.
\end{align*}

\newpage

\section{Tables for the optimal ratio in various polytopes}
The following tables give the optimal ratio together with $\delta_\parallel$ for regular polytopes in 2, 3, 4, and 5 dimensions with \begin{align*}
    \phi = \frac{\sqrt{5}+1}{2}.
\end{align*} The values of $r_{\text{opt}}$ and $\delta_\parallel$ have been experimentally verified using the Overlap Testing Python programs provided in \cite{programs}. 

\renewcommand{\arraystretch}{1.5}

\begin{table}[h!]
\centering
\begin{tabular}{|c|c|c|}
\hline
\textbf{Polytope (2D)} & \(\mathbf{\delta_\parallel}\) & \(\mathbf{r_{\text{opt}}}\) \\
\hline
Triangle             & $\ell$ & $\frac{1}{2}$ \\
Square               & $\ell$ & $\frac{1}{2}$ \\
Pentagon             & $\phi \ell$ & $\frac{1}{\phi}$ \\
Hexagon              & $2\ell$ & $\frac{2}{3}$ \\
\hline
\end{tabular}
\caption{Optimal Ratio for Regular 2D Polytopes}
\end{table}

\begin{table}[h!]
\centering
\begin{tabular}{|c|c|c|}
\hline
\textbf{Polytope (3D)} & \(\mathbf{\delta_\parallel}\) & \(\mathbf{r_{\text{opt}}}\) \\
\hline
Tetrahedron            & $\ell$ & $\frac{1}{2}$ \\
Cube                   & $\ell$ & $\frac{1}{2}$ \\
Octahedron             & $\ell$ & $\frac{1}{2}$ \\
Icosahedron           & $\phi \ell$ & $\frac{1}{\phi}$ \\
Dodecahedron           & $(\phi + 1)\ell$ & $\frac{\phi + 1}{\phi + 2}$ \\
\hline
\end{tabular}
\caption{Optimal Ratio for Regular 3D Polytopes}
\end{table}

\begin{table}[h!]
\centering
\begin{tabular}{|c|c|c|}
\hline
\textbf{Polytope (4D)} & \(\mathbf{\delta_\parallel}\) & \(\mathbf{r_{\text{opt}}}\) \\
\hline
5-cell (4-simplex)      & $\ell$ & $\frac{1}{2}$ \\
8-cell (4-cube)         & $\ell$ & $\frac{1}{2}$ \\
16-cell                & $\ell$ & $\frac{1}{2}$ \\
24-cell                & $2\ell$ & $\frac{2}{3}$ \\
\hline
\end{tabular}
\caption{Optimal Ratio for Regular 4D Polytopes}
\end{table}

\begin{table}[h!]
\centering
\begin{tabular}{|c|c|c|}
\hline
\textbf{Polytope (5D)} & \(\mathbf{\delta_\parallel}\) & \(\mathbf{r_{\text{opt}}}\) \\
\hline
5-simplex              & $\ell$ & $\frac{1}{2}$ \\
5-cube                 & $\ell$ & $\frac{1}{2}$ \\
5-orthoplex            & $\ell$ & $\frac{1}{2}$ \\
\hline
\end{tabular}
\caption{Optimal Ratio for Regular 5D Polytopes}
\end{table}

\newpage

\section{Applications}
\begin{figure}[h!]
    \centering
    \includegraphics[width=0.5\linewidth]{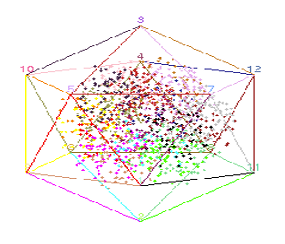}
    \caption{GCG played on a regular icosahedron with \(r = 0.5\) \cite{annie}}
    \label{fig:annie}
\end{figure}

The applications of this paper are primarily in the field of biology. In \cite{annie}, Thomas discusses a three-dimensional Chaos Game representation of protein sequences. In the paper, a GCG is played on a regular icosahedron using a ratio of \(r = 0.5\). The resulting graph (Figure \ref{fig:annie}) reveals no visible structure.

\begin{figure}[h!]
    \centering
    \includegraphics[width=0.5\linewidth]{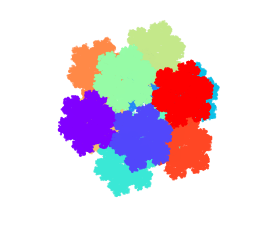}
    \caption{GCG played on a regular icosahedron with \(r = \frac{1}{\phi}\) \cite{programs}}
    \label{fig:icosahedron}
\end{figure}

However, if Thomas had used the optimal ratio of \(r = \frac{1}{\phi}\), a discernible structure would have emerged, resembling the one shown in Figure \ref{fig:icosahedron}.

Similarly, in \cite{SUN20201904}, Sun et al. play a GCG on a regular dodecahedron using \(r = \frac{\sqrt{5}-1}{\sqrt{5}+2\sqrt{3}-1} \approx 0.7370\). The authors correctly state that this ratio ensures the scaled-down dodecahedra do not intersect. However, for the dodecahedra to also touch, the ratio should be \(r = \frac{\phi + 1}{\phi + 2} \approx 0.7236\).

\newpage

\printbibliography

\end{document}